\documentclass[11pt]{amsart}

\usepackage{hyperref}

\begin{document}
\bibliographystyle{plain}

\newtheorem{conjecture}{Conjecture}
\newtheorem{thm}{Theorem}[section]
\newtheorem{lemma}[thm]{Lemma}

\title{Algorithms to Uniformly Generate Factored Smooth Integers}

\author{Eric Bach}
\address{Computer Sciences Department \\
  University of Wisconsin-Madison}
\email{bach@cs.wisc.edu} 
\author{Jonathan Sorenson}
\address{Computer Science \&\ Software Engineering Department \\
  Butler University} 
  \email{jsorenso@butler.edu} 
\date{\today\\Poster presented at ANTSXIV, 2020}

\maketitle
\begin{abstract}
Let $x\ge y>0$ be integers.
A positive integer is $y$-smooth if all its prime divisors are
  at most $y$.
Let $\Psi(x,y)$ count the number of $y$-smooth integers up to $x$.
We present several algorithms that will generate an integer $n\le x$
  at random, with known prime factorization, such that
  $n$ is $y$-smooth.
We begin by describing algorithms to compute $\Psi(x,y)$ exactly and
  to enumerate $y$-smooth integers up to $x$ in
  lexicographic order by prime divisor.
Both of these are based on Buchstab's identity, and were likely
  known before.
Then we present
an algorithm that accepts as input a parameter $r$, $0\le r<1$,
  and returns the integer $n$ that is at position $\lfloor r\Psi(x,y)\rfloor$ 
  in the lexicographic ordering of all $y$-smooth integers up to $x$.
Here position 0 is the first position.
Thus, $n$ is generated uniformly so long as $r$ is chosen uniformly.
This algorithm has a running time of $O(\Psi(x,y)\log\log y)$ arithmetic
  operations.
We then explore the tradeoff between speed and rigor.
By relaxing the uniformity of the output and allowing for
  multiple heuristics in our runtime analysis, 
  we improve the running time to
$$
  O\left( \frac{ (\log x)^3 }{\log\log x} \right)
$$
arithmetic operations.
We conclude with a sample run by generating a $10000$-smooth
  integer $\le 10^{100}$.

\end{abstract}

\nocite{FGKT2016,dj71,brent,suzuki2004,suzuki2006,vdLW69}
\nocite{Alladi87,Hildebrand87,Hensley87,EK40}

\section{Introduction}

Given integers $x$ and $y$, with $x\ge y>0$,
  we say a positive integer $n\le x$ is $y$-\textit{smooth} if every prime
  divisor of $n$ is $\le y$.
(In the literature, \textit{friable} is sometimes used in place of
  \textit{smooth}.)

In this paper, we describe and analyze new algorithms to
  generate a positive integer $n\le x$ that is $y$-smooth,
  together with a complete list of $n$'s prime divisors.
Further, we want $n$ to be chosen uniformly at random 
  from among all $y$-smooth integers $\le x$.
Let $\Psi(x,y)$ count the number of $y$-smooth integers $\le x$.
Then we want the probability a particular $n$ is chosen to be $1/\Psi(x,y)$.

To achieve this, we start by describing an algorithm that lists all
  $y$-smooth integers $\le x$.
This method is recursive,
  and uses Buchstab's identity as its central control mechanism.
We then make use of a random input $r$, with $0\le r<1$,
  and return the smooth number at position $k=\lfloor r\Psi(x,y)\rfloor$
  from that listing.
If we view the algorithm's execution as a tree, 
  then each leaf is a smooth number,
  and the listing is composed of all the leaves.
We then prune the algorithm so that its execution generates the leaf for 
  the smooth number at position $k$, 
  while doing as little extra work as possible.

We then analyze the running time of our algorithm,
  and look at tradeoffs we can make to improve performance.
In particular, we look at relaxing the uniformity requirement,
  we consider assuming unproven conjectures like the ERH
  and Cram\'er's conjecture on the maximum size of gaps between primes,
  and we look at using randomized prime tests.
We also discuss average-case versus worst-case running times.


In particular, we show the following.
\begin{enumerate}
\item
Given integers $x\ge y>0$ and $r$ chosen uniformly at random with $0<r\le 1$,
there is an algorithm to produce the $y$-smooth integer $n$, along with
  $n$'s complete prime factorization,
  that occurs at position $\lfloor r\Psi(x,y)\rfloor$ from the lexicographic
  enumeration of all $y$-smooth integers up to $x$.
This algorithm has a running time of $O(\Psi(x,y)\log\log y)$ arithmetic
  operations.
This is Theorem \ref{thm:slow}.
\item
With the same inputs $x,y,r$ as above,
there is an algorithm to produce a $y$-smooth $n$ along with its factorization,
  where $n$ occurs at position $r\Psi(x,y)(1+o(1))$ in the lexicographic
  enumeration,
  so our uniformity condition has been relaxed.
This algorithm has a running time that is roughly linear in $y$.
This is Theorem \ref{thm:withprimes}.
\item
With the same inputs $x,y,r$ as above,
under several heuristic assumptions which we specify later,
there is an algorithm to produce a $y$-smooth $n$ along with its factorization,
  where $n$ occurs at position $r\Psi(x,y)(1+o(1))$ in the lexicographic
  enumeration.
This algorithm has a running time of 
$O( (\log x)^3/\log\log x )$ arithmetic operations.
This is Theorem \ref{thm:fast}.
\end{enumerate}

\subsection{Model of Computation\label{sec:model}}

We measure algorithm complexity by counting the number
  of arithmetic operations on integers of $O(\log x)$ bits.
This includes addition/subtraction, multiplication/division with remainder,
  and other basic operations such as array indexing, comparisons,
  and branching.
As a consequence, for example, performing a base-$2$ strong pseudoprime test
  on an integer $\le x$ would take $O(\log x)$ arithmetic operations.

We also work with floating-point real numbers with up to $O(\log x)$ bits
  of precision with exponents bounded by $O(\log x)$,
  and again each operation on such a number takes unit time.

In practice, we would expect integers $\le y$ to fit into 64-bit machine
  words, and floating-point numbers would be represented using
  a $\texttt{long double}$ data type, which is typically 64 or 80 bits.
Integers between $y$ and $x$ are not used often, and so
  multi-precision software (such as GMP) should not be needed except in
  rare cases.
In fact, it might be possible to represent $x$ in a double precision floating
  point representation.

\subsection{Paper Outline}
The rest of our paper is organized as follows:
We begin in \S\ref{sec:enum} with an algorithm to compute $\Psi(x,y)$
  exactly, and then modify it to enumerate all $y$-smooth integers $\le x$.
Then in \S\ref{sec:slow} we show how to modify
  the algorithm to selectively list just one $y$-smooth integer $\le x$
  by pruning our first algorithm.
In \S\ref{sec:tradeoff} we look for ways to improve the runtime performance,
  including estimating values of $\Psi$.
In \S\ref{sec:anal} we give a more detailed running time analysis of
  our algorithm, making use of heuristics and relaxing our uniformity
  requirement.
We then discuss some special cases and applications in \S\ref{sec:special},
  and conclude with an example run in \S\ref{sec:run}.
Our code for the sample run is available here:
  \url{https://github.com/sorenson64/rsi}.

\section{Enumerating All Smooth Integers\label{sec:enum}}

\subsection{Buchstab's Identity}
As mentioned in the Introduction, our algorithms are based on
  a version of Buchstab's identity,
\begin{equation} \label{eq:buch}
\Psi(x,y) = 1 + \sum_{p\le y} \Psi(x/p,p) ,
\end{equation}
which decomposes $\Psi(x,y)$ by its largest prime divisor
  \cite[\S5.3]{Tenenbaum}.
Combine this with some base cases,
\begin{enumerate}
  \item $\Psi(x,y) = 0$ if $x<1$ or $y<1$, and
  \item $\Psi(x,1) = 1$ if $x\ge 1$, and
  \item $\Psi(x,2) = \lfloor \log_2 x \rfloor +1$ if $x\ge 1$, and
  \item $\Psi(x,y) = \lfloor x\rfloor$ if $y\ge x$, and
  \item $\Psi(x,y) = \Psi( \lfloor x \rfloor,\lfloor y \rfloor )$,
\end{enumerate}
and you have a simple recursive algorithm to compute the exact
  value of $\Psi(x,y)$.
We have the following.

\begin{lemma}\label{lemma:psi}
Given integers $x\ge y>0$,
there is an algorithm to compute the exact value of $\Psi(x,y)$
using $O(\Psi(x,y))$ arithmetic operations.
\end{lemma}
\begin{proof}
The only piece missing from the discusson above is that we can
  find the primes up to $y$ using the Atkin-Bernstein sieve in
  $O(y/\log\log y)$ arithmetic operations \cite{AB2004},
  but $\Psi(x,y)\ge \lfloor y\rfloor $ when $x\ge y$,
  making this negligible.
\end{proof}

\subsection{Example}
We illustrate this using $\Psi(15,3)$.
We have
$$
  \Psi(15,3)=1+\Psi(15/2,2)+\Psi(15/3,3).
$$
Here $\Psi(15/2,2)=\Psi(7,2)=\lfloor \log_2 7 \rfloor+1=3$,
leaving $\Psi(15/3,3)=\Psi(5,3)$ as a recursive case:
$$
 \Psi(5,3)=1+\Psi(5/2,2)+\Psi(5/3,3)=1+(\lfloor \log_2 2.5 \rfloor+1) + 1
   = 4.
$$
Plugging back in, we get $\Psi(15,3)=8$, which is correct.

To have this method generate a list of numbers, we merely need to
  remember the prime $p$ from the sum in Buchstab's identity
  used to construct the recursive call.
We will annotate these primes as subscripts to $\Psi$ and redo our example.
We have
$$
  \Psi(15,3)=1+\Psi_2(15/2,2)+\Psi_3(15/3,3).
$$
This generates the $1$.
The call to $\Psi(7,2)$ generates powers of $2$, 
  namely $2^0,2^1,2^2$.
The $2$ coming in via the subscript tells us multiply
  everything generated by $2$, giving us the list
  $2\cdot 2^0,2\cdot 2^1,2\cdot 2^2$.
$\Psi(5,3)$ would recursively 
  generate the list $1,2,3,4=2^2$, and applying the $3$
  subscript gives the list $3,3\cdot2, 3\cdot3, 3\cdot2^2$.
The complete list, then, is
$$ (1, 2\cdot 2^0,2\cdot 2^1,2\cdot 2^2, 3,3\cdot2, 3^2, 3\cdot2^2), $$
or $1,2,4,8,3,6,9,12$.

\subsection{Intuition}
If we take Buchstab's identity and divide through by $\Psi(x,y)$,
  we see that when constructing a number $n\le x$,
  each prime $p\le y$ has a chance to be $n$'s largest prime divisor,
  and that probability is $\Psi(x/p,p)/\Psi(x,y)$.
No prime is chosen, or $n=1$, with probability $1/\Psi(x,y)$.

Also note that if we omit the 4th base case 
  ($\Psi(x,y)=\lfloor x\rfloor$ if $y\ge x$),
  the enumeration produced is in lexicographic order based 
  on decreasing prime factorization.
So from our example, the enumeration would be written like this:
\begin{quote}
  $1$ \\
  $2$ \\
  $2 \cdot 2$ \\
  $2 \cdot 2 \cdot 2$ \\
  $3$ \\
  $3 \cdot 2$ \\
  $3 \cdot 2 \cdot 2$ \\
  $3 \cdot 3$
\end{quote}
Here's another example starting at position $100$ in the enumeration
  of $7$-smooth integers $\le 1000$:
\begin{quote}
$7\cdot 3\cdot 3  = 63$ \\
$7\cdot 3\cdot 3\cdot 2  = 126$ \\
$7\cdot 3\cdot 3\cdot 2\cdot 2  = 252$ \\
$7\cdot 3\cdot 3\cdot 2\cdot 2\cdot 2  = 504$ \\
$7\cdot 3\cdot 3\cdot 3  = 189$ \\
$7\cdot 3\cdot 3\cdot 3\cdot 2  = 378$ \\
$7\cdot 3\cdot 3\cdot 3\cdot 2\cdot 2  = 756$ \\
$7\cdot 3\cdot 3\cdot 3\cdot 3  = 567$ \\
$7\cdot 5  = 35$ \\
$7\cdot 5\cdot 2  = 70$ \\
$7\cdot 5\cdot 2\cdot 2  = 140$ \\
\end{quote}
If one is curious, $\Psi(1000,7)=141$.

\subsection{An Enumeration Algorithm}
Next, we describe an algorithm to enumerate $y$-smooth integers.
The stack $S$ is used to store the primes we wrote in the example 
  as subscripts to $\Psi$,
  and the vector/array $V$ holds the enumeration, so
  $V=[1,2,4,8,3,6,9,12]$ from our example, or, with complete factorizations,
  $$V=[1,[2],[2,2],[2,2,2],[3],[3,2],[3,3],[3,2,2]].$$

\begin{tabbing}MM\=MM\=MM\=MM\=\kill
Procedure \textbf{Enumerate}($x$,$y$,$S$,$V$): \+\\[5pt]
  If $x<1$ Then do nothing and return; \\[3pt]
  If $S$ is empty Then Append 1 onto $V$; \\
  Else Append a copy of $S$ onto $V$; \\[3pt]
  For each prime $p\le y$ Do:\+\\
    Push $p$ onto $S$; \\
    \textbf{Enumerate}($x/p$,$p$,$S$,$V$); \\
    Pop $p$ from $S$;
\end{tabbing}
Our example, then, would be produced using \textbf{Enumerate}($15,3,S=[],V=[]$).
The $S$ and $V$ arguments must be pass-by-reference, 
  but $x$ and $y$ should be ordinary pass-by-value arguments.

\begin{lemma}\label{lemma:enum}
Procedure \textbf{Enumerate}() will list all $y$-smooth integers up to $x$,
along with a list of their prime divisors, with multiplicity,
in $O(\Psi(x,y)(\log \log x + \frac{\log x}{\log y})$ arithmetic operations.
\end{lemma}
\begin{proof}
It should be clear that the running time is proportional to
  the total number of primes stored in the vector $V$ at the end.
The cost of finding the primes up to $y$ is negligible, as
  discussed earlier.
Due to the work of Alladi \cite{Alladi87}, Hensley \cite{Hensley87},
  and Hildebrand \cite{Hildebrand87} we know the average number
  of prime divisors is
\begin{equation} \label{eq:divcnt}
 O\left(\log\log x + \frac{\log x}{\log y} \right).
\end{equation}
So our running time is $O( \Psi(x,y) (\log\log x+ (\log x)/\log y ))$
  arithmetic operations.
\end{proof}

To generate one randomly chosen smooth number, we simply construct the
  list $V$ and choose one entry uniformly at random.

\subsubsection*{Algorithm 1}
\begin{enumerate}
  \item Set $S:=[]$ and $V:=[]$.
  \item Call \textbf{Enumerate}($x,y,S,V$).
  \item Choose $r$ uniformly at random, with $0\le r < 1$.
  \item Set $k:=\lfloor r\Psi(x,y) \rfloor 
               = \lfloor r\times \mbox{length}(V) \rfloor$.
  \item Output $V[k]$.
\end{enumerate}

This gives us the following theorem.

\begin{thm}\label{thm:alg1}
Given integers $x\ge y>0$ and a real parameter $r$ with $0\le r < 1$,
Algorithm 1 returns a $y$-smooth integer $n\le x$ in factored form
  from position $k=\lfloor r \Psi(x,y)\rfloor$ in the 
  lexicographic enumeration of all $y$-smooth integers $\le x$,
and its running time is
  $O( \Psi(x,y) (\log\log x+ (\log x)/\log y ))$
  arithmetic operations.
\end{thm}

So long as $r$ is selected uniformaly at random, then so is
  the output $n$ from Algorithm 1.

Note that the list $V$ requires $O(\Psi(x,y)\log x)$ bits of space.

\section{Pruning Algorithm 1\label{sec:slow}}

Algorithm 1 above works, but is quite slow because it constructs
  all smooth numbers before selecting one at random,
  and as noted above, uses a large amount of space.
We only need to construct one.
So, to improve our algorithm, 
  instead of working through all the primes $p\le y$
  in the main loop of the \textbf{Enumerate} function,
  we use the value of $k=\lfloor r\Psi(x,y) \rfloor$ 
  to figure out the prime $p\le y$ on which to recurse.
From Buchstab's identity, this means finding the 
  consecutive primes $p_1$ and $p_2$
  such that
  $$
    1+\sum_{p\le p_1} \Psi(x/p,p) < k+1 \le 1+\sum_{p\le p_2} \Psi(x/p,p),
  $$
  or, rewriting this applying Buchstab's identity,
  $$
    \Psi(x,p_1) < k+1 \le \Psi(x,p_2).
  $$
  Then, we recurse on the branch that computes $\Psi(x/p_2,p_2)$.
  The value of $k$ for the recursive call, 
    $k^\prime$, is simply $k-\Psi(x,p_1)$.
The understanding, here, is that if $k=0$ then we are returning $1$,
  but if $k+1\le \Psi(x,2)$, then we recurse on $p_2=2$, which means
  in this special case, $p_1$ is set to $1$, even though $1$ is not prime.

This gives the following structure for our algorithm.
\subsubsection*{Algorithm 2}
\begin{enumerate}
  \item Set $S:=[]$.
  \item Choose $r$ uniformly at random, with $0\le r < 1$.
  \item Set $k:=\lfloor r\Psi(x,y) \rfloor$.
  \item Output \textbf{Branch}($x,y,k,S$).
\end{enumerate}
\begin{tabbing}MM\=MM\=MM\=MM\=\kill
Function \textbf{Branch}($x$,$y$,$k$,$S$): \+\\[3pt]
  If $x<1$ Then return nothing/null; \\[3pt]
  If $k=0$ Then \+\\
    If $S$ is empty Then return 1; \\
    Else return $S$; \-\\[3pt]
  Find consecutive primes $p_1<p_2$ such that
    $\Psi(x,p_1)<k+1\le \Psi(x,p_2)$. \\[3pt]
  Push $p_2$ onto $S$; \\
  Set $k^\prime=k-\Psi(x,p_1)$; \\
  Return \textbf{Branch}($x/p_2$,$p_2$,$k^\prime$,$S$); 
\end{tabbing}
The running time of Algorithm 2 is at worst
  proportional to the recursion depth of the \textbf{Branch} function,
  times the time for one execution of the \textbf{Branch} function
  (excluding the recursive call at the end).
The recursion depth is exactly the number of prime divisors, with
  multiplicity, in the number $n$ generated by the algorithm.
From (\ref{eq:divcnt}) above, this is $O(\log\log x + (\log x)/\log y)$
  in the average case, and $O(\log x)$ in the worst case.
It should be easy to see that the bottleneck in one \textbf{Branch} execution 
  is finding the primes $p_1,p_2$, the \textit{Find} step,
  which we will discuss in detail below.
Note that the computation of
  $k^\prime=k-\Psi(x,p_1)$ can re-use the value of $\Psi(x,p_1)$ from
  the Find step.

If we use our exact algorithm to compute the value of $\Psi(x,y)$ and
  find all primes up to $y$ using a sieve as before, we obtain the following.
\begin{thm}\label{thm:slow}
Given integers $x\ge y>0$, Algorithm 2 will output a $y$-smooth integer $n\le x$
  with its complete prime factorization, and $n$ will be chosen uniformly
  at random from among all $y$-smooth integers $\le x$.
Algorithm 2 takes $O(\Psi(x,y)\log\log y)$ arithmetic operations.
\end{thm}
\begin{proof}
With a list of primes up to $y$ in hand, interpolation search can find
  $p_1,p_2$ in $O(\log\log y)$ steps.
It remains to show that all evaluations of the $\Psi$ function take time
  $O(\Psi(x,y)\log\log y)$ where $x,y$ are the original inputs to the algorithm.
Write $n=q_1q_2\cdots q_m$ where the $q_i$ are the prime divisors of $n$ 
  with $q_i\ge q_{i+1}$ for $1\le i<m$.
At the top level of the algorithm $q_1$ is chosen to recurse on,
  and it takes $O(\Psi(x,y)\log\log y)$ time to do the evaluations of $\Psi$.
In the recursive call, $x$ is replaced by $x/q_1$.
At this next level, the total time spent evaluating $\Psi$ is
  $O(\Psi(x/q_1,y)\log\log y)$.
At the level below that,
  $O(\Psi(x/(q_1q_2),y)\log\log y)$.
From Theorem 12 in \cite[\S5.5]{Tenenbaum}, we have
  $$\frac{1}{q}\Psi(x,y) \le \Psi(x/q,y)(1+o(1)),$$
  giving us
$$
\sum_{i=1}^m \Psi(x/(q_1\cdots q_i),y) 
  \le \Psi(x,y)(1+o(1)) \sum_{i=1}^m \frac{1}{q_1\cdots q_i}
  = O(\Psi(x,y)).
$$
This completes the proof.
\end{proof}
Although this is not much of an improvement over Algorithm 1 
  in terms of running time,
  its use of space is quite a bit better with no need for the
  vector $V$.
It needs only $O(y+\log x)$ bits, mostly to store the primes up to $y$.


\section{Trading Rigor for Speed}\label{sec:tradeoff}

The algorithms we have presented so far are quite slow.
There are two ways to go faster:
\begin{enumerate}
\item
If we are willing to compromise on the uniformity condition of the output,
  we can use approximations for $\Psi$ instead of exact values.
Given how slow computing exact values for $\Psi$ is, this should have
  a noticeable impact.
  (See \S\ref{sec:psi}.)
\item
If we hope to get a running time significantly below linear in $y$,
  we do not have the time to find all the primes up to $y$.
Instead, we must somehow find the primes $p_1,p_2$ as needed. 
To do this, we break this into two steps: Find a value of $t$
  such that the estimate $\Psi(x,t)-k$ is near zero, and then
  find consecutive primes $p_1<t\le p_2$.
  (See \S\ref{sec:search} for finding $t$ and
  \S\ref{sec:primes} for finding $p_1,p_2$ from $t$.)
\end{enumerate}

\subsection{Algorithms to Estimate $\Psi(x,y)$\label{sec:psi}}

In addition to our recursive method above for exactly computing $\Psi(x,y)$,
  there are a number of algorithms to estimate $\Psi$ in the literature:
\begin{itemize}
  \item A method to compute $\Psi(x,y)$ exactly that is potentially
    faster than using Buchstab's identity
    \cite{Bernstein95-2},
  \item Methods to give upper and lower bounds on $\Psi(x,y)$ using
    formal power series 
    \cite{Bernstein2002,PS06},
  \item Methods based on the saddle-point approach of Hildebrand and Tenenbaum 
    \cite{HS97,LP2018,Sorenson2000,suzuki2004,suzuki2006}, 
    (see also \cite{BS13} for a Lagarias-Miller-Odlyzko-style algorithm),
    and
  \item Methods based on the Dickman-deBruijn function
    \cite{Dickman30,vdLW69}.
  \nocite{HT93}
\end{itemize}
We can imagine circumstances in which each of these would be useful.
It turns out that currently, the fastest method is the oldest,
\begin{equation} \label{eq:xrho}
  \Psi(x,y) \approx x\cdot \rho(u),
\end{equation}
where $\rho(u)$ is the Dickman-deBruijn function, 
  and $u=u(x,y)=(\log x)/\log y$.
To be precise,
Hildebrand proved that
\begin{equation} \label{eq:hildebrand}
 \Psi(x,y) = x \rho(u) \left( 1+ O(\frac{\log(u+1)}{\log y}) \right) 
\end{equation}
under the condition that $\log y > (\log \log x)^{5/3+\epsilon}$
\cite[Cor.\ 9.3, \S5.5]{Tenenbaum}.
Under the assumption of the ERH, 
  the range on $y$ can be extended down to $y>(\log x)^{2+\epsilon}$.

We will use this estimate so long as $y$ is in the correct range,
  because it asymptotically estimates $\Psi$ correctly,
  and the $\rho$ function is fast and easy to compute.

\subsubsection{Computing $\rho$ for Large $y$}
The idea is to pre-compute a table of values $\rho(u)$
  for all $u$ up to a bound $B$ we will determine later
  (although we do know $B\le \log_2 x$),
  using the trapezoid rule (see \cite{vdLW69}).

We know that the absolute error on a subinterval of size $h$
  is bounded by $(1/12)\rho^{\prime\prime}(u_h)\cdot h^3$,
  where $u_h$ is a point in the interval of size $h$
  (See \cite[\S7.2]{CdB}).
From \cite[Cor.\ 8.3,\S5.4]{Tenenbaum} we have that
  $\rho^{\prime\prime}(u) \sim (\log u)^2 \rho(u)$,
  giving us a relative error of $(1+ O(h^2/(\log u)^2))$.
If we choose $h=1/\log x$, then the relative error on our
  computation of $\rho$ will be smaller than the relative error
  in Hildebrand's result above.
For details, see \cite{vdLW69,MZM89,BP96}.

This will take $O( B (\log x)^2 )$ arithmetic operations to build the table,
  which will contain $B\log x$ floating point numbers.
Computing estimates for $\Psi(x,y)$ will then involve only a table
  lookup and a multiplicaton, or constant time.

\subsubsection{The Cutoff $L(x)$}

Define $L(x)$ to be the cutoff point where, if $y<L(x)$,
  then the $x\rho(u)$ estimate in (\ref{eq:hildebrand}) is no longer valid.
So then
\begin{eqnarray*}
L(x)&=&(\log x)^{2+\epsilon} \mbox{\ if we assume the ERH, and} \\
L(x)&=&\exp\left[(\log \log x)^{5/3+\epsilon}\right] \mbox{\ otherwise.}
\end{eqnarray*}
If $u\le (\log x)/\log L(x)$, then (\ref{eq:hildebrand}) is valid,
  and we use (\ref{eq:xrho}) to estimate $\Psi(x,y)$.
This also means we can take $B=B(x)=(\log x)/\log L(x)$ 
  when computing our table of $\rho$ values.
Note that as the algorithm progresses, $x$ gets smaller, and $B(x)$
  decreases as $x$ does.
Also, observe that when a value of $u$ arises at any point in our computation,
  if $u>B$, then $(\ref{eq:hildebrand})$ is not valid, and we don't
  need a value of $\rho$ for this $u$.

\subsubsection{Saddle-Point Methods for Small $y$}

We are all set for when $y\ge L(x)$.
When $y<L(x)$, we will use a theorem
  of Hildebrand and Tenenbaum based on the saddle-point method to
  estimate $\Psi$ \cite[Thm 10,\S5.5]{Tenenbaum}.
This allows $y$ to be as small as $2$.

If we are assuming the ERH,
  then we will use the HT-fast algorithm of \cite{Sorenson2000},
  which gives a running time of 
  $O( (\sqrt{y}/\log y) \log\log x)$ assuming we have a list of
  primes up to $\sqrt{y}$.
Without the ERH,
  we use algorithm HT from \cite{HS97} with 
  Lagarias-Miller-Odlyzko (LMO) summation as described
  in \cite{BS13} for a running time of $y^{2/3+o(1)}\log\log x$.
If we have a list of primes up to $y$ available,
  the non-LMO version has a running time of
  $O( y  (\log\log x)/\log y )$.

\subsubsection{When $x,y$ are Both Small}

If we discover that $\Psi(x,y)\log x$ is smaller than the time already spent
  constructing $n$,
  we can revert to using Algorithm 1 to finish off the computation.
When $x,y$ are both small, the asymptotics break down and $n$ is much less
  likely to be chosen suitably uniformly, and so the exact Algorithm 1 is
  preferable in this case.


\subsection{Searching for $p_1,p_2$}\label{sec:search}

Next, we look at the \textit{Find} step.
Our first task is to address the special case when $p_1=1,p_2=2$.
This only takes constant time, since $\Psi(x,2)=\lfloor \log_2 x \rfloor+1$.
So going forward, we can assume $p_1>1$.

Since we are using an approximation algorithm to compute $\Psi$,
  and we may not have a list of primes up to $y$ available,
  as stated above,
  we break this down into two steps:
\begin{enumerate}
\item Find a real number $t$ such that $\Psi(x,t)-(k+1)$ is near zero.
\item Find consecutive primes $p_1<p_2$ with $p_1<t\le p_2$.
\end{enumerate}

\subsubsection{Finding $t$ when $t\ge L(x)$}

We can compute $\Psi(x,L(x))$ using the $x\rho(u)$ method in constant time
  to quickly determine whether our solution $t$ 
  is larger or smaller than $L(x)$.
Our first case is when $t\ge L(x)$.
Writing $v=v(x,t)=(\log x)/\log t$, this means
 we can use the $x\rho(v)$ estimate for the rest of this step.

The estimate $x\rho(v)$ is continuous, strictly increasing in $t$,
  and differentiable so long as $v\ge1$.
Thus, bisection or binary search works well, and since
  evaluations of $\Psi$ take constant time, the overall cost
  is $O(\log y)$ arithmetic operations. 
We will see this is negligible
  compared to other steps in the algorithm.

That said, we can use Newton's method to
  find $t$ in at most $O(\log\log y)$ evaluations of $\Psi$.
We only need to know $t$ to the nearest integer.
To be precise, we are finding a root $t$ to the function 
  $f(t)=\Psi(x,t)-(k+1)$.
Using the defining equation $\rho^\prime(v)=-\rho(v-1)/v$,
  we obtain the iteration function 
$$
  g(t) = t-\frac{f(t)}{f^\prime(t)}= t-\frac{\rho(v)-(k+1)/x}{\rho(v-1)} t\log t.
$$
If $v$ is large enough, it is easy to prove quadratic convergence.
See \cite[\S3.5]{CdB}.

\subsubsection{Finding $t$ when $t<L(x)$}

The estimates based on the saddle-point method involve
  sums over primes, and so will not be continuous in general.
However, $f(t)=\Psi(x,t)-(k+1)$ 
  is still a non-decreasing function of $t$ with a simple root,
  so we can use the Illinois algorithm \cite{dj71}
  or Brent's algorithm \cite{brent},
  both of which converge super-linearly,
  to find a root in $O(\log\log y)$ evaluations of $\Psi$.
Once the size of the interval to search has shrunk
  to length $2(\log y)^2$, we switch to bisection to avoid
  any issues with discontinuity, and with no effect on the
  asymptotic running time.
See also
  \cite[\S3.3,\S3.5]{CdB}.

\subsubsection{Computing $p_1,p_2$ from $t$\label{sec:primetests}}

With $t$ in hand, we want to find the largest prime $p_1<t$ and the
  smallest prime $p_2\ge t$.
If we happen to have a list of primes up to and a bit past $t$,
  we can do a quick interpolation search on the list to find $p_1,p_2$
  in $O(\log\log t)$ time, which is negligible compared to finding $t$ itself.
So for now we will assume we don't have access to such a list.

Prime testing can be expensive, so to minimize the number of prime tests,
  we do the following.
Set an interval length $w=2\lceil \log t \rceil$,
  and sieve an interval of length $w$ with midpoint $t$
  by the primes $\le \log t$.
We then perform a base-2 strong pseudoprime test on the
  $O(w/\log\log t)$ integers that pass the sieve.
This takes $O(w (\log t)/\log\log t)$ time.
Identify candidate values for $p_1,p_2$, if they exist on the interval,
  and prime test them.
If two such candidates are found that pass, then we are done.
If not, set $w:=2\cdot w$ and repeat the process
  until two candidates are found.
Of course, if one candidate is found but not the other,
  just continue doubling the interval above or below $t$, but not both,
  as appropriate.

Non-prime integers that pass the base-2 strong pseudoprime
  test are very rare; see \cite{PSW80} for example.
So on average, only two prime tests are needed.
Also, by the prime number theorem, we expect to find these primes
  without having to double $w$ at all, or maybe just a constant number of times.
So on average, the cost is $O( (\log t)^2/\log\log t)$ plus
  a constant number of prime tests.

In the worst case, every base-2 strong pseudoprime test requires a follow-up
  full prime test.
However,
  we cannot guarantee a prime will show up until $w\gg{t}^{0.525}$
  \cite{BHP2001}.
If we are sieving an interval of length $\sqrt{t}$ or larger, then we
  may as well sieve by primes up to $\sqrt{t}$ so that no prime tests
  are required.
So in the worst case, the cost is $O({t}^{0.525}/\log\log t)$ to run
  the Atkin-Bernstein sieve on an interval of size $O({t}^{0.525})$.


For a more reasonable worst-case running time,
  we can choose to use Cram\'er's conjecture.
\begin{conjecture}[Cram\'er]\label{conj1}
  There exist absolute constants $n_0,c>0$ such that if $n>n_0$,
  then $p_{n+1}-p_n \le c ( \log p_n )^2$.
\end{conjecture}
See \cite{FGKT2016} for a discussion 
  of this conjecture and other related results on gaps between primes.
With this conjecture, we have $w=O( (\log t)^2 )$
  for a bound of $O( (\log t)^2/\log\log t)$ prime tests.

\subsection{Prime Testing\label{sec:primes}}

In practice, it makes sense to use a fast, 
  probabilistic prime test such as Miller-Rabin \cite{Miller76,Rabin80},
  to find $p_1,p_2$ and then optionally verify their primality using
  a more rigorous test.
Since the Miller-Rabin test fails with probability at most $1/4$,
  we can use each test at most $O(\log\log x)$ times, on each of $p_1,p_2$,
  and still get the error
  probability down to $o(1)$ over the entire algorithm.
A single test is $O(\log y)$ operations to perform, the same as
  a base-2 strong pseudoprime test, giving a cost of
  $O( (\log y) \log\log x)$ operations for each candidate prime.

Possible options for the more rigorous prime test include
  the following:
\begin{itemize}
  \item If we are assuming the ERH, Miller's test \cite{Miller76}
    takes $O( (\log y)^3 )$ arithmetic operations.
  \item If we allow random numbers, Bernstein's variant of the AKS test
     \cite{Bernstein2004} takes expected $(\log y)^{3+o(1)}$
     arithmetic operations.
    Note that this is much slower in practice than Miller's algorithm.
    Also note that the randomness is only in the running time, not correctness.
  \item If we don't allow the ERH or random number use,
    but have access to a sufficiently large table of pseudosquares,
    the pseudosquares prime test \cite{LPW96}
    is fast at $O((\log y)^2)$ arithmetic operations.
    (See \cite{Sorenson06,Sorenson10a} for algorithms to build the required
     table.)
  \item
    And then there's the AKS test \cite{AKS04}
    which has a variant due to Lenstra and Pomerance that takes
    $(\log y)^{5+o(1)}$ arithmetic operations.
    No ERH or random numbers are required.
\end{itemize}

If we are in a theoretical situation where random bits are a scarce resource,
  building a table of pseudosquares is probably the way to go.

\section{Analysis and Tradeoffs\label{sec:anal}}

The running time of our algorithm is proportional to the following
  expression:
$$
  T + D\cdot( S\log\log y + P )
$$
where
\begin{itemize}
\item $T$ is the cost to build the table of $\rho(u)$ values,
\item $D$ is the number of prime divisors or the recursion depth,
\item $S$ is the cost of evaluating $\Psi(x,y)$, and
\item $P$ is the cost of finding $p_1,p_2$ from $t$.
\end{itemize}
Let's look at each one in detail.

\subsection{Building a table of $\rho(u)$ values}

We know $T=O( (\log x)^3/\log L(x) )$, and $L(x)$ depends only on
  whether we assume the ERH or not:
\begin{eqnarray*}
  T &=& O\left(\frac{(\log x)^3}{\log\log x}\right) \quad \mbox{(with ERH)} \\
  T &=& O\left(\frac{(\log x)^3}{(\log\log x)^{5/3+\epsilon}}\right)
      \quad \mbox{(without ERH)} 
\end{eqnarray*}
Note that if we are generating large numbers of random smooth factored
  integers, computing the $\rho$ table can be viewed as one-time preprocessing.

\subsection{Counting prime divisors}

From (\ref{eq:divcnt}) we have
\begin{eqnarray*}
  D &=& O(\log x) \quad \mbox{(worst case),} \\
  D &=& O\left(\log\log x + \frac{\log x}{\log y}\right)
      \quad \mbox{(average case).} 
\end{eqnarray*}

\subsection{Evaluating $\Psi(x,y)$}

As we look at $S$, we'll find it helpful to define $z$ as the
  smaller of $y$ and $L(x)$.
When $y\ge L(x)$, we know $S$ is constant time.
When $y<L(x)$, if we are assuming the ERH, then
  we use algorithm HT-fast, and we know that $L(x)=(\log x)^{2+\epsilon}$,
  so that
  $$S=O\left(\sqrt{z}\cdot\frac{\log\log x}{\log z}\right)
         \le (\log x)^{1+\epsilon+o(1)}
    \quad\mbox{(with ERH)}.$$
Note that we could drop the $o(1)$ in the exponent here by choosing
  a different $\epsilon$, but we leave it there to remind the
  reader that we are masking lower-order multiplicative terms.
Without the ERH, we use the LMO version of algorithm HT, and
   $L(x)=\exp[ (\log\log x)^{5/3+\epsilon}]$,
  giving 
 $$S=O( z^{2/3+o(1)}\log\log x ) \le  
    \exp\left[ (\log\log x)^{5/3+\epsilon+o(1)}\right]
    \quad\mbox{(without ERH)}. $$
Again, we leave the $o(1)$ to indicate the masking of lower order factors.
If we happen to have a list of primes up to $y$ available, 
we may as well use Algorithm HT, giving
 $$S=O( (z/\log z)\log\log x ) \le  
    \exp\left[ (\log\log x)^{5/3+\epsilon+o(1)}\right]
    \quad\mbox{(without ERH)}. $$

\subsection{Finding $p_1,p_2$ from $t$}

Our last one, $P$, is the most interesting.

If we use an average-case analysis,
  then $P=O( (\log y)^2/\log\log y + M )$, where the first term is
  the cost of sieving and using base-2 strong pseudoprime tests, and $M$
  is the cost of one prime test.
If we are assuming the ERH, we can set $M=O((\log y)^3)$, or better yet,
  if we are using probabilistic prime tests, $M=O( (\log y)\log\log x)$.
With neither the ERH nor randomization, we use the AKS test
  for $M=(\log y)^{5+o(1)}$:
\begin{eqnarray*}
P&=& O\left( \frac{ (\log y)^2 }{\log\log y} \right)
  \quad\mbox{(avg.case, random)}, \\
P&=& O\left( (\log y)^3 \right)
  \quad\mbox{(avg.case, ERH, not random)}, \\
P&=& (\log y)^{5+o(1)}
  \quad\mbox{(avg.case, no ERH, not random)}.
\end{eqnarray*}

If we use a worst-case analysis, we have more options to consider.
If we assume Conjecture \ref{conj1}, then
  $P=O(M (\log y)^2/\log\log y)$.
If we further allow the use of random numbers, 
    we can use a probabilistic prime test.
In this case, the expected (not average!) 
    running time is $O(\log y)$ to reject a number as not prime.
This gives $P=O( (\log y)^3/\log\log y )$, with the final two probabilistic
  prime tests taking only $O( (\log y) \log\log x)$ time.
Without randomization, $M=O((\log y)^3)$ with the ERH,
  giving $P=O( (\log y)^6/\log\log y))$,
  and $M=(\log y)^{5+o(1)}$ without the ERH, for
  $P= (\log y)^{7+o(1)}$:
\begin{eqnarray*}
P&=& O\left( \frac{(\log y)^3}{\log\log y} \right)
  \quad\mbox{(worst case, Conj.\ref{conj1}, random)}, \\
P&=& O\left( \frac{ (\log y)^6 }{\log\log y} \right)
  \quad\mbox{(worst case, Conj.\ref{conj1}, ERH, not random)}, \\
P&=& (\log y)^{8+o(1)}
  \quad\mbox{(worst case, Conj.\ref{conj1}, no ERH, not random)}.
\end{eqnarray*}

If we don't assume Conjecture \ref{conj1} but do assume the RH,
  then we will sieve an interval of size $O(\sqrt{y}\log y)$ for primes,
  searching for $t$ on that interval quickly finds $p_1,p_2$, 
  giving us 
$$P=O\left(\frac{\sqrt{y}(\log y)}{\log\log y}\right) 
  \quad\mbox{(worst case, no Conj.\ref{conj1}, RH)}.
$$
With neither conjecture, the approach is the same but the length of
  the interval is slightly larger, giving
$$P= O\left( \frac{y^{0.525}}{\log\log y} \right)
  \quad\mbox{(worst case, no Conj.\ref{conj1}, no RH)},
$$
  from \cite{BHP2001}.

\subsection{Complexity Results}
Given all our options for conjectures, the use of randomness,
  and worst-case versus average-case, we can obtain a wide range of
  running times for our algorithm.
Rather than go through all of these, we present three
  versions that seem reasonably useful.

This first theorem follows from choosing to find all primes $\le y$
  and then using Algorithm HT to estimate $\Psi$.
If a list of primes is already available, the first $y/\log\log y$ term
  can be dropped from the running time.
\begin{thm}\label{thm:withprimes}
  Let $x\ge y>0$ be integers, and let $r\in[0,1)$ be chosen
    uniformly at random.
  Algorithm 2 will construct a integer $n$ with known prime factorization such
    that $n\le x$ and $n$ possesses no prime divisor larger than $y$.
  The worst-case running time of this algorithm is
  $$
   O\left( \frac{y}{\log\log y} +
     \frac{y}{\log y} (\log\log x)(\log\log y) 
     \left(\frac{\log x}{\log y} + \log\log x \right)
   \right)
  $$
  arithmetic operations.

  Here $n$ is chosen asymptotically uniformly, in the sense that as
    $x,y\rightarrow\infty$, $n$'s relative position in the full 
    lexicographic enumeration
    of $y$-smooth integers $\le x$ is
    $\lfloor r\Psi(x,y)(1+o(1))\rfloor$.
\end{thm}

This second theorem follows from choosing our options 
  so as to obtain the smallest
  possible running time, namely ERH, randomized prime tests,
  and an average-case running time analysis.
\begin{thm}\label{thm:fast}
  Let $x\ge y>0$ be integers, and let $r\in[0,1)$ be chosen
    uniformly at random.
  Under the assumption of the Extended Riemann Hypothesis (ERH),
    Algorithm 2 will construct a integer $n$ with known factorization such
    that $n\le x$ and $n$ possesses no prime divisor larger than $y$.
  The divisors of $n$ are all prime with probability $1-o(1)$.
  The average running time of this algorithm is
  $$
    O\left( \frac{(\log x)^3}{\log\log x} \right)
  $$
  arithmetic operations.
  If a table of values of $\rho(u)$ is precomputed for $u$ up to
    $O( (\log x)/\log\log x)$,
  then the running time drops to 
  $$O\left( \frac{(\log x)(\log y)}{\log\log y}
    + \frac{ (\log x)^{2+\epsilon}}{\log y} \right) $$
  arithmetic operations, where $\epsilon>0$.

  As in the sense of Theorem \ref{thm:withprimes}, $n$ is chosen
    asymptotically uniformly.
\end{thm}

Our third theorem is as follows.
\begin{thm}\label{thm:safe}
  Let $x\ge y>0$ be integers, and let $r\in[0,1)$ be chosen
    uniformly at random.
  Algorithm 2 will construct a integer $n$ with known factorization such
    that $n\le x$ and $n$ possesses no prime divisor larger than $y$.
  As in the sense of Theorem \ref{thm:withprimes}, $n$ is chosen
    asymptotically uniformly.
  Let $\epsilon>0$ and let $L(x)=\exp[ (\log\log x)^{5/3+\epsilon}]$
    as required by (\ref{eq:hildebrand}).
  Let $z$ be the smaller of $y$ and $L(x)$.
  Under the assumption of Conjecture \ref{conj1},
    Algorithm 2 has a worst-case running time of
  $$
    z^{2/3+o(1)}(\log x)\log\log x 
       + O\left( \frac{(\log x)^3}{ (\log\log x)^{5/3+\epsilon}}\right).
  $$
\end{thm}
Here the proof follows from our discussion above, with the no-ERH option,
  a worst-case running time, with Conjecture \ref{conj1}, and
  no randomized prime tests (they would not help the overall running time).
We are using the LMO version of Algorithm HT to estimate $\Psi$
  when $y<L(x)$.
The second term is from precomputing $\rho(u)$.

\section{Special Cases\label{sec:special}}

\subsection{Bach's Result\label{sec:bach}}
We can use Algorithm 2 to produce random factored numbers 
  somewhat similarly to 
  Bach's algorithm \cite{Bach88} by simply setting $y=x$.
The running time would be $O((\log x)^3/\log\log x)$ to precompute the
  table of $\rho$ values, and then $O((\log x)^2/\log\log x)$ operations
  on average to generate each integer $n$ (using Theorem \ref{thm:fast}).
Bach's algorithm uses $O(\log x)$ prime tests.
From section \ref{sec:primetests}, Bach's algorithm would give an
  $O((\log x)^2\log\log x)$ running time if we use probabilistic
  prime tests, so our new method is faster by a factor 
  proportional to $(\log\log x)^2$ if we ignore precomputation costs.
Note that this generates integers in $[1,x]$ instead of $(x/2,x]$,
  and Bach's algorithm guarantees an exact uniform distribution, 
  whereas ours does not.

\subsection{Semismooth Numbers\label{sec:semi}}
The central control structure of our algorithm is Buchstab's identity.
With some straightforward modifications,
  we could adapt our approach to make use
  of the generalized Buchstab identity given in \cite{BS13}
  for generating random factored semismooth integers.

\section{Example Run\label{sec:run}}

We implemented our algorithm in C++ on a linux desktop workstation
  and ran it with $x=10^{100}$, $y=10000$, and $r=0.5$.
It generated the following list of prime divisors for $n$:
\begin{quote}
2 3 5 7 29 31 97 113 113 113 157 223 241 503 509 569 691 727 1033 1367 1571 2141 2339 2617 2741 3041 3221 3547 3989 4021 4513 4999 5573 6577 7573 9463
\end{quote}
The resulting $n$ is roughly $4.29\cdot 10^{97}$, which occupies
  a position near $2.05\cdot 10^{61}$ in the enumeration.
The run took less than $0.35$ seconds of wall time.

Since this $y$ is small enough, we found all primes $\le y$
  and used Algorithm HT to estimate $\Psi(x,y)$,
  or in other words, the version of Theorem \ref{thm:withprimes}.

Our code can be found here:
\url{https://github.com/sorenson64/rsi}.

\bibliography{all}
\end{document}